\documentclass[12pt]{iopart}
\usepackage{amsfonts}
\usepackage{amssymb}

\providecommand{\U}[1]{\protect\rule{.1in}{.1in}}
\providecommand{\U}[1]{\protect\rule{.1in}{.1in}}

\newtheorem{definition}{Definition}

\newtheorem{proposition}{Proposition}
\newtheorem{remark}{Remark}

\newenvironment{proof}[1][Proof]{\noindent\textbf{#1.} }{\ \rule{0.5em}{0.5em}}

\begin{document}

\title {Nonlocal interpretation of $\lambda$-variational symmetry-reduction method}

\author{D Catalano Ferraioli$^1$, P Morando$^2$}

\address{$^1$Dipartimento di Matematica, Universit\`{a} di Milano, via Saldini 50, I-20133 Milano, Italy}
\address{$^2$Istituto di Ingegneria Agraria, Facolt\`{a} di Agraria, Universit\`{a} di Milano, Via Celoria, 2 - 20133 Milano, Italy}

\ead{catalano@mat.unimi.it and paola.morando@polito.it}

\begin{abstract}

In this paper we give a geometric interpretation of a reduction method based on the so called $\lambda$-variational symmetry (C. Muriel, J.L. Romero and P. Olver 2006 \emph{Variational $C^{\infty}$-symmetries and Euler-Lagrange equations} J. Differential equations \textbf{222} 164-184). In general this allows only a partial reduction but it is particularly suitable for the reduction of variational ODEs with a lack of computable local symmetries. We show that this method is better understood as a nonlocal symmetry-reduction.
\end{abstract}

\pacs{02.30.Hq} \ams{37K05, 34C14, 70S10}
\vspace{2pc} \noindent{\it Keywords}:  $\lambda$-variational symmetry, local
symmetries, nonlocal symmetries, ODE
reduction, variational symmetries\\


\section{Introduction}
In the last few decades various generalizations of the notion of classical symmetry of ODEs have been proposed. It is well known, in fact,  that finding all local symmetries for an ODE is not always possible and one may encounter equations solvable by quadratures but with a lack of local symmetries. Hence local symmetries are sometimes inadequate and various attempts, for a more effective symmetry-reduction method, have been done. Among the proposed generalizations, in our opinion, a special attention is deserved by those introduced in \cite{MuRo,MRO06}.\\
The method introduced in \cite{MuRo} is based on the notion of $\lambda$-symmetry. In fact, if an equation is invariant under a $\lambda$-symmetry, one can obtain a complete set of functionally independent invariants and reduce the order of the equation by one as for Lie symmetries. Despite their name, however, $\lambda$-symmetries are neither Lie point nor higher symmetries; nevertheless, as shown in \cite{Ca}, $\lambda$-symmetries of an ODE $\mathcal{E}$ can be interpreted as shadows of some nonlocal symmetries. This interpretation of $\lambda$-symmetries has the main advantage that it allows to reinterpret the $\lambda$-symmetry reduction as a particular case of the standard symmetry reduction method. A detailed discussion on $\lambda$-symmetries and further generalized methods can be found in \cite{Ca,GaMo,Paola07,ultimo}.\\
In this paper we are mainly concerned with the geometric interpretation of a reduction method based on an alternative notion of variational symmetry, which was introduced in \cite{MRO06} and is known as $\lambda$-variational symmetry. In particular we use a nonlocal framework to interpret any variational $\lambda$-symmetry $X$ as a standard variational symmetry $\tilde{X}$ for an associated variational problem. As a consequence, we can show that to any $X$ there corresponds a standard first integral $\tilde{I}$ for the associated variational problem. However $\tilde{I}$ depends on a nonlocal variable and, in order to use it in the reduction of the initial Euler-Lagrange equations $\mathcal{E}$, one has to restrict $\mathcal{E}$ on the submanifold $\{\tilde{I}=0\}$. This is a kind of conditional reduction for the Euler-Lagrange equations.\\
Above interpretation clarifies the geometric aspects of the method given in \cite{MRO06} and provides explicit invariant formulas for the first integral $\tilde{I}$. We do not discuss computability issues neither special applications of above method, hence the example at the end of last section is just a simple illustration of the above reduction procedure.
\section{Preliminaries}
We assume that the reader is familiar with the geometry of differential
equations (see for example \cite{Vin-et-al,Olv1}) and we only collect
basic facts and notations we use in the paper. Let $B$ and $M$ be
smooth manifolds and $\pi:M\rightarrow B$ be a $q$-dimensional bundle.
We denote by $\pi_{k}:J^{k}(\pi)\rightarrow B$ the $k$-order \emph{jet
bundle} associated to $\pi$. Since we are only concerned with the
case $\dim B=1$, we assume that $B$ and $M$ have local coordinates
$x$ and $(x,u^{1},...,u^{q})$, respectively. Moreover, we only consider
second order systems of Euler-Lagrange equations for regular first
order Lagrangians. In the standard jet coordinates $(x,u^{a},u_{1}^{a},u_{2}^{a})$
on $J^{2}(\pi)$, such a second order system of ODEs can be written
in the following form \begin{equation}
\mathcal{E}:=\left\{ \begin{array}{l}
u_{2}^{a}=f^{a}(x,u,u_{1}): \;\; a=1,...,q \end{array}\right\}.\label{equazione}\end{equation}

Geometrically $\mathcal{E}$ is interpreted as a submanifold of $J^{2}(\pi)$
and is naturally equipped with the contact distribution generated
by the restriction $\bar{D}_{x}$ of the total derivative $D_{x}$
to $\mathcal{E}$:
\begin{equation}
\bar{D}_{x}=\partial_{x}+u_{1}^{a}\partial_{u^{a}}+f^{a}\partial_{u_{1}^{a}}.\label{D_x_ristr}\end{equation}
 Here and throughout the paper we use Einstein summation convention
over repeated indices.

The infinite prolongation of (\ref{equazione}) corresponds to the
following submanifold of the infinite jet bundle $J^{\infty}(\pi)$
\[
\mathcal{E}^{\infty}:=\left\{ \begin{array}{l}
D_{x}^{i}\left(u_{2}^{a}-f^{a}(x,u,u_{1})\right)=0:\;a=1,...,q\;\;\mathrm{and} \;\;i=0,1,...\end{array}\right\}.\]

 We denote by $\mathcal{C}$ and $\bar{\mathcal{C}}$ the \emph{contact distributions} on $J^{\infty}(\pi)$ and $\mathcal{E}^{\infty}$, respectively. Notice that both $\mathcal{E}^{\infty}$ and $\mathcal{E}$ have
the same finite dimension and the contact distribution on $\mathcal{E}^{\infty}$
is $1$-dimensional and still generated by the vector field $\bar{D}_{x}$.
The contact distribution on $\mathcal{E}^{\infty}$ can also be described
as the annihilator space of the contact ideal $\bar{\mathcal{C}}$
generated by the $1$-forms $\{du^{a}-u_{1}^{a}dx,du_{1}^{a}-f^{a}dx:a=1,...,q\}$.

In the paper we only consider
the so called \emph{external symmetries} of $\mathcal{E}$, i.e., symmetries
of the form \begin{equation}
X=\xi\partial_{x}+\eta_{i}^{a}\partial_{u_{i}^{a}}\label{**}\end{equation}
 with $\xi,\eta^{a}$ arbitrary smooth function on $\mathcal{E}^{\infty}$
and\[
\eta_{i}^{a}=D_{x}(\eta_{i-1}^{a})-u_{i}^{a}D_{x}(\xi),\qquad\eta_{0}^{a}=\eta^{a}.\]

Hence, $X$ is the infinite prolongation of a \emph{Lie point (or contact)
symmetry} iff $\xi,\eta_{0}^{a}$ are functions on $M$ (or $J^{1}(\pi)$,
respectively).

On $J^{\infty}(\pi)$, symmetries $X$ of $\mathcal{C}$ which are
tangent to $\mathcal{E}^{(\infty)}$ are called \emph{higher symmetries}
of $\mathcal{E}$ and are determined by the condition $\left.X(F)\right\vert _{\mathcal{E}^{(\infty)}}=0$.

When working on infinite jets spaces, since $D_{x}$ is a trivial
symmetry of $\mathcal{C}$, it is convenient to gauge out from higher
symmetries the terms proportional to $D_{x}$. This leads to consider
only symmetries in the so called \emph{evolutive form}, i.e., symmetries
of the form $X=D_{x}^{i}(\varphi^{a})\partial_{u_{i}^{a}}$, $\varphi^{a}:=\eta^{a}-u_{1}^{a}\xi$.
The functions $\varphi^{a}$ are called the \emph{generating functions}
(or characteristics) of $X$.

\bigskip

A $1$-dimensional \emph{covering} for (\ref{equazione}) is the infinite
prolongation $\tilde{\mathcal{E}}:=(\mathcal{E}^{\prime})^{\infty}$
of a system of the form\begin{equation}
\mathcal{E}^{\prime}:\left\{ \begin{array}{l}
u_{2}^{a}=f^{a}(x,u,u_{1}),\\
w_{1}=\lambda,\end{array}\right.\label{covering}\end{equation}
 with $\lambda$ a smooth function on $\mathcal{E}^{(\infty)}$.
Equation $\mathcal{E}^{\prime}$ can be interpreted as a submanifold
of the second jet space of a trivial bundle $\tau:\mathbb{R}^{2+q}\rightarrow\mathbb{R}$
with local coordinates $(x,u^{a},w)$. The contact distribution on
$\tilde{\mathcal{E}}$ is the $1$-dimensional distribution generated
by the vector field \begin{equation}
\tilde{D}_{x}=\bar{D}_{x}+\lambda\partial_{w}.\label{D_x_tilde}\end{equation}
 The contact distribution on $\tilde{\mathcal{E}}$ can also be described
as the annihilator space of the contact ideal $\tilde{\mathcal{C}}$
generated by the $1$-forms $dw-\lambda dx$ and $\{du^{a}-u_{1}^{a}dx,du_{1}^{a}-f^{a}dx:a=1,...,q\}$.

\emph{Nonlocal symmetries} of $\mathcal{E}$ are symmetries of the vector
field (\ref{D_x_tilde}) and can be determined through a symmetry
analysis of the system $\left(\mathcal{E}^{\prime}\right)^{(\infty)}$
on $J^{\infty}(\tau)$. Therefore nonlocal symmetries
of $\mathcal{E}$ have the form \begin{equation}
Y=\xi\partial_{x}+\eta_{i}^{a}\partial_{u_{i}^{a}}+\psi_{i}\partial_{w_{i}}\label{pippo}\end{equation}
 with $\eta^a_{i}=\widetilde{D}_{x}(\eta^a_{i-1})-\widetilde{D}_{x}(\xi)u^a_{i}$
and $\psi_{i}=\widetilde{D}_{x}(\psi_{i-1})-w_{i}\widetilde{D}_{x}(\xi)$.

\medskip

An interesting example of nonlocal symmetry occurring in literature
is related to the notion of \emph{$\lambda$-symmetry} for an ODE $\mathcal{E}$
(see \cite{MuRo}).

If $\lambda$ is a smooth function on $J^{1}(\pi)$,
then we say that the $\lambda$-prolongation to $J^{k}(\pi)$ of a vector field
$X=\xi\partial_{x}+\eta^{a}\partial_{u^{a}}$ on $M$ is
the vector field $X^{[\lambda,k]}=\xi\partial_{x}+\eta_{\lbrack\lambda,i]}^{a}\partial_{u_{i}^{a}}$
with \[
\eta_{\lbrack\lambda,0]}^{a}=\eta^{a},\qquad\eta_{\lbrack\lambda,i]}^{a}=D_{x}(\eta_{\lbrack\lambda,i-1]}^{a})-D_{x}(\xi)u_{i}^{a}+\lambda\left(\eta_{\lbrack\lambda,i-1]}^{a}-\xi u_{i}^{a}\right).\]
 Moreover, we say that a vector field $X^{[\lambda,k]}$ is a $\lambda$-symmetry of $\mathcal{E}$ if and only
if $X^{[\lambda,k]}$ is tangent to $\mathcal{E}$.

Despite their name, $\lambda$-symmetries are neither Lie
symmetries nor higher symmetries of $\mathcal{E}$. Nevertheless,
as discussed in \cite{Ca}, $\lambda$-symmetries can be interpreted
as shadows of nonlocal symmetries. More precisely, $\mathcal{E}$
admits a $\lambda$-symmetry $X$ iff $\mathcal{E}^{\prime}=\{u_{k}^{a}=f^{a},w_{1}=\lambda\}$,
with $\lambda\in C^{\infty}(\mathcal{E})$, admits a (higher) symmetry
with generating functions of the form $\varphi^{\alpha}=e^{w}\varphi_{0}^{\alpha}(x,u,u_{1},...,u_{k-1})$,
$\alpha=1,...,q+1$.

Since $\lambda$-symmetries are of great interest in the applications,
and analogously their nonlocal counterparts, it is convenient to introduce
the notion of \emph{$\lambda$-covering}: if a covering system $\mathcal{E}^{\prime}$,
defined by (\ref{covering}),
admits a nonlocal symmetry $Y$ with generating functions \begin{equation}
\varphi^{\alpha}=e^{w}\varphi_{0}^{\alpha}(x,u,u_{1},...,u_{k-1}),\qquad\alpha=1,..,q+1\label{gen_func_lambda}\end{equation}
 then $\mathcal{E}^{\prime}$ will be called a $\lambda$-covering
for $\mathcal{E}$ defined by (\ref{equazione}).

\bigskip{}

\section{Main results}

The relation between standard symmetries of a Lagrangian and conservation
laws for the corresponding Euler-Lagrange equations is described by classical Noether theorem.

An extension of this classical result, to the case of $\lambda$-variational symmetries, was proposed in \cite{MRO06}. In this section we give an interpretation of this result in terms of nonlocal symmetries and we point out that to any $\lambda$-variational symmetry can be associated a standard variational symmetry (with corresponding standard conservation law) for the singular Lagrangian $\tilde{L}=\tau^*(L)$ in a $\lambda$-covering of the Euler-Lagrange equations.

If $L:J^{1}(M)\rightarrow\mathbf{{R}}$ is a first order regular Lagrangian,
we denote by \begin{equation}
\Theta\ =\ \frac{\partial L}{\partial u_{1}^{a}}\ (du^{a}-u_{1}^{a}dx)\ +\ L\ dx\label{PC_form_1ord}\end{equation}
 the Poincaré-Cartan 1-form associated to $L$.

\begin{definition} \medskip{}
\ If (\ref{equazione}) is an Euler-Lagrange equations with a regular
first order Lagrangian $L(x,u,u_{1})$, a vector field \ $X$ is called a (divergence) $\lambda$-variational symmetry for
$L$ iff \[
X(L)+L(D_{x}+\lambda)\xi=(D_{x}+\lambda)R,\]
 with $X$ the $\lambda$-prolongation of a vector field $X_{0}=\xi \partial_x +\eta^a \partial_{u^a}$
on $M$, and $R$ a function on $M$. \end{definition}

\begin{remark} Notice that, in general, a $\lambda$-variational symmetry
is neither a $\lambda$-symmetry nor a symmetry for the corresponding
Euler-Lagrange equations. \end{remark}

One can show (see \cite{Paola07}) that

\begin{proposition} A $\lambda$-prolonged vector
field $X$ is a (divergence) $\lambda$-variational symmetry for a Lagrangian $L(x,u,u_{1})$
iff \begin{equation}
L_{X}(\Theta)+(X\lrcorner\Theta-R)\lambda dx-dR\in\bar{\mathcal{C}}.\label{Prop_3_1}\end{equation}
\end{proposition}

If $\tau$ denotes the projection of a $\lambda$-covering for the Euler-Lagrange equations $\mathcal{E}$ of a Lagrangian $L$, we write $\tilde{L}=\tau^*(L)$ and with a slight abuse of notation  $\Theta=\tau^*(\Theta)$. The function $\tilde{L}$ is a singular Lagrangian on $J^1(\tau)$.

Now we prove that to any $\lambda$-variational symmetry it can be
associated a standard variational symmetry in the covering space $\tilde{\mathcal{E}}$.

\begin{proposition} The vector field $X$ is a (divergence) $\lambda$-variational symmetry
for $L$ iff in the covering (\ref{covering}), for any smooth function $A$ of $(x,u,u_{1},w)$, the vector field $\tilde{X}=e^{w}(X+A\partial_{w})$
is such that

\begin{equation}
L_{\tilde{X}}(\Theta)-d(e^{w}R)\in\tilde{\mathcal{C}}.\end{equation}
\end{proposition}

\begin{proof} In view of the form of $\tilde{X}$, one has \begin{equation}
L_{\tilde{X}}\Theta-d(e^{w}R)=e^{w}L_{X+A\partial_{w}}\Theta+d(e^{w})\wedge(X+A\partial_{w})\lrcorner\Theta-Rde^{w}-e^{w}dR.\label{Prop_4_1}\end{equation}
 On the other hand, since $\Theta$ do not depends on $w$, \[
L_{X+A\partial_{w}}\Theta=L_{X}\Theta,\qquad(X+A\partial_{w})\lrcorner\Theta=X\lrcorner\Theta.\]
 Hence \[
\begin{array}{l}
L_{\tilde{X}}\Theta-d(e^{w}R)=e^{w}[L_{X}\Theta+(X\lrcorner\Theta)dw-Rdw-dR]=\\
e^{w}[L_{X}\Theta+(X\lrcorner\Theta-R)\lambda dx-dR]+e^{w}(X\lrcorner\Theta-R)(dw-\lambda dx),\end{array}\]
 and then one readily gets the thesis in view of (\ref{Prop_3_1})
and the fact that $dw-\lambda dx\in\tilde{\mathcal{C}}$.\\
 \end{proof}

Notice that, with respect to the singular Lagrangian $\tilde{L}$, $\tilde{X}$ is a standard variational symmetry in the covering.
Nevertheless $\tilde{X}$ is not a symmetry for the corresponding Euler-Lagrange equations, in view
of the singularity of the Lagrangian function.

\medskip

It is well known that Noether theorem associates to any standard (divergence)
variational symmetry $X$ for $L$ (with divergence term $R$) the
first integral $I=X\lrcorner\Theta-R$ for the corresponding Euler-Lagrange
equations (\ref{equazione}). Unfortunately, the same is not true
for (divergence) $\lambda$-variational symmetries. In fact, if $X$ is a $\lambda$-variational symmetry for $L$, $I=X\lrcorner\Theta-R$ is not a first
integral for (\ref{equazione}), but satiesfies (see \cite{MRO06,Paola07,GCvar})

\begin{equation}
\bar{D}_{x}(I)+\lambda I=0\hspace{0.2in}mod\;\mathcal{E}.\label{***}\end{equation}

Nevertheless, in the covering system (\ref{covering}) one has the
following

\begin{proposition} \label{Prop_7}Given a $\lambda$-variational symmetry
$X$ for $L$, then $\tilde{I}=\tilde{X}\lrcorner\Theta-e^w R=e^w I$
is a first integral for $\tilde{D}_{x}$ (i.e., for the covering system
(\ref{covering})). \end{proposition}

\begin{proof} Since $X$ is a $\lambda$-variational symmetry, $I=X\lrcorner\Theta-R$
satisfies (\ref{***}). Then the thesis follows by the fact that $\tilde{D}_{x}(X\lrcorner\Theta-R)=\bar{D}_{x}(X\lrcorner\Theta-R)$
and the following direct computation:\[
\begin{array}{l}
\tilde{D}_{x}(\tilde{X}\lrcorner\Theta-e^{w}R)=\tilde{D}_{x}(e^{w}(X+A\partial_{w})\lrcorner\Theta-e^{w}R)\\
\qquad\hspace{0.31in}=\tilde{D}_{x}(e^{w})X\lrcorner\Theta+e^{w}\tilde{D}_{x}(X\lrcorner\Theta)-\lambda e^{w}R-e^{w}\tilde{D}_{x}(R)\\
\qquad\hspace{0.31in}=e^{w}\left(\lambda X\lrcorner\Theta+\bar{D}_{x}(X\lrcorner\Theta)-\lambda R-\bar{D}_{x}(R)\right)=e^{w}(\bar{D}_{x}(I)+\lambda I)\end{array}\]
 then $\tilde{D}_{x}(\tilde{I})$ vanishes in view of (\ref{***}).
\end{proof}

Since $\tilde{I}$ is a first integral of $\tilde{\mathcal{E}}$,
it follows that one can use $\tilde{I}=c$ to reduce by one the order
of $\tilde{\mathcal{E}}$. On the contrary, in view of (\ref{***}),
$I$ is not a first integral of $\mathcal{E}$, and one only has that

\begin{equation}
\bar{D}_{x}(I)=0\;\;\;\; mod\;\mathcal{E}\cup\{I=0\}.\end{equation}
 Therefore equation $\mathcal{E}$ can be reduced only on the hypersurface
$\{I=0\}$. This partial reduction of $\mathcal{E}$ is a conditional
reduction and coincide with that proposed in the paper \cite{MRO06}.

\bigskip

\noindent\textbf{Example.} Equation $\mathcal{E}$ defined by

\begin{equation}
u_2=\frac{x^{2}\ln u}{u}-\frac{x^{2}}{u}+\ln u-1\label{eq:example}\end{equation}
is the Euler-Lagrange equations corresponding to the Lagrangian function

\begin{equation}
L=\frac{{\it u_{1}}^{2}}{2}+{\it xu_{1}}(1-\ln u)+x^{2}\ln u \left(\frac{\ln u}{2}-1\right).
\end{equation}

\noindent Above equation does not admit Lie point symmetries, hence it cannot
be reduced by standard symmetry reduction techniques. Nevertheless,
it admits $X=\partial_{u}$ as a $\lambda$-variational symmetry with
$\lambda=x/u$. Then one can show that the first integral $\tilde{I}$
of the $\lambda$-covering system $\tilde{\mathcal{E}}$ reads\[
\tilde{I}=e^{w}\left(u_{1}-x\ln u+x\right).\]
Then the reduced ODE of (\ref{eq:example}) is \begin{equation}
u_{1}-x\ln u+x=0.\label{eq:red_example}\end{equation}
One can readily show that any solution of (\ref{eq:red_example})
is also a solution of (\ref{eq:example}).

\bigskip{}

\end{document}